\documentclass{article}
\usepackage[left=2cm,top=2cm,right=3cm,bottom=1cm, footskip=.5cm]{geometry}
\usepackage{amsmath}
\usepackage{graphicx}
\usepackage{tikz}
\usepackage{fancyvrb}
\usepackage{pdfpages}
\usepackage[makeroom]{cancel}
\usepackage{amssymb} 
\usepackage[parfill]{parskip}
\usepackage{mathtools}
\usepackage{hyperref}
\usepackage{tikz-cd}
\usepackage{setspace}
\usepackage{amsthm}
\newtheoremstyle{colon}%
{}
{}
{\itshape}
{}
{\bfseries}
{:}
{ }
{}

\theoremstyle{colon}

\newtheorem{Lemma}{Lemma}[subsection]
\newtheorem{Theorem}[Lemma]{Theorem}

\newtheorem{Corollary}[Lemma]{Corollary}
\newtheorem{Definition}[Lemma]{Definition}

\title{A double coset formula for the genus of a nilpotent group}

\begin{document}
	
	\date{}
	\author{A. Ronan}
	\maketitle
	
\begin{abstract}
We derive double coset formulae for the genus and extended genus of a finitely generated nilpotent group G, using the notions of bounded and bounded above automorphisms of $\prod G_S$, which are defined relative to a fixed fracture square for G.
\end{abstract}
	
\section{A double-coset formula for the genus of a nilpotent group}

\subsection{Introduction}

Let $T,S$ and, for each $i$ in some indexing set $I$, $T_i$ be sets of primes such that $T = \cup_i T_i$ and $T_i \cap T_j = S$ for all $i \neq j$. Suppose also that $T \neq S$. Throughout we let $G$ be an $f\mathbb{Z}_T$-nilpotent group and consider a fixed diagram:

\[ \begin{tikzcd}
G \arrow{r}{(\psi_i)} \arrow[swap]{d}{\sigma} & \prod G_{T_i} \arrow{d}{\phi} \arrow{dr}{\prod \phi_i} \\
G_S \arrow[swap]{r}{\omega} & (\prod G_{T_i})_S \arrow[swap]{r}{\tilde{\pi}} & \prod G_S \\
\end{tikzcd}
\]

where each $\psi_i$ is a localisation at $T_i$, $\phi$ is a localisation at $S$, $\sigma$ is a localisation at $S$, $\phi_i$ is the unique localisation at $S$ such that $\phi_i \psi_i = \sigma$, $\omega$ is the localisation of $(\psi_i)$ and $\tilde{\pi}$ is the unique map making the triangle on the right commute. It follows from these definitions that $\tilde{\pi} \omega = \Delta$. In $[1, 7.5]$ a map was defined which sends an automorphism $\alpha \in \prod Aut(G_S)$ to the pullback of $\alpha \circ \prod \phi_i$ along $\Delta$. It was claimed that this map was a surjection onto the extended genus of $G$, where the extended genus of $G$ is the set of isomorphism classes of T-local nilpotent groups $H$ such that $H_{T_i} \cong G_{T_i}$ for every $i \in I$. However, it is not necessarily true that the localisations of the pullback group at each $T_i$ will agree with $G_{T_i}$ as the following example shows:

\[ \begin{tikzcd}
\mathbb{Z} \arrow{r}{(\psi_i)} \arrow[swap]{d}{\sigma} & \prod \mathbb{Z}_{\{p_i\}}  \arrow{d}{\prod \phi_i} \\
\mathbb{Q} \arrow[swap]{r}{\Delta} & \prod \mathbb{Q} \\
\end{tikzcd}
\]  

where each of the undefined maps is the inclusion sending $1$ to $1$. Now consider the automorphism $\alpha = \prod p_i$ of $\prod \mathbb{Q}$, where $p_i$ denotes multiplication by the ith prime number $p_i$. Then the image of $\alpha \circ \prod \phi_i$ consists of elements $(q_i) \in \prod \mathbb{Q}$ such that if $q_i = \frac{a_i}{b_i}$ with $a_i,b_i$ coprime, then $a_i$ is divisible by $p_i$. In particular, the image of $\alpha \circ \prod \phi_i$ intersects the image of $\Delta$ only at $0$. Therefore, the pullback group of $\alpha \circ \prod \phi_i$ along $\Delta$ is $0$ which does not localise to $\mathbb{Z}_{\{p_i\}}$ for each $i$. 

Nevertheless, this example turns out to be instructive. Suppose, instead, that $\alpha = \prod(\frac{u_i}{v_i})$ with $u_i$ and $v_i$ coprime integers. Suppose that $\alpha$ is 'bounded' in the sense that there are only finitely many primes which divide some $u_i$ or $v_i$. Then the induced pullback is isomorphic to $\mathbb{Z}$, which is the unique abelian group in the genus of $\mathbb{Z}$, where the genus is the subset of the extended genus consisting of finitely $T$-generated groups. If, instead, $\alpha$ is only 'bounded above' in the sense that there are only finitely many primes which divide some $u_i$, then the pullback turns out to be in the extended genus of $\mathbb{Z}$ and the maps $\psi_i$ are localisations at $T_i$. In fact, the pullback group will not be finitely generated unless $\alpha$ is 'bounded'. Note that in the counterexample we formulated, the map $\alpha$ was neither 'bounded' nor 'bounded above'. 

With this in mind, the purpose of these notes will be to prove the following pair of double coset results relating to the genus and extended genus of $G$ respectively:

\textit{\textbf{Theorem:} The genus of $G$ is in 1-1 correspondence with the double coset:}
	
	\begin{center}
		$Aut(G_S) \text{\textbackslash} Aut_{b}(\prod_i G_S) / \prod_i Aut(G_{T_i})$
	\end{center} 
	
\textit{where $Aut_{b}(\prod_i G_S)$ is the subgroup of automorphisms of the form $\prod_i \alpha_i$ which are $S$-bounded, see Definition 1.2.2}

\textit{\textbf{Theorem:} The extended genus of $G$ is in 1-1 correspondence with the double coset:}
	
	\begin{center}
		$Aut(G_S) \text{\textbackslash} Aut_{b.a.}(\prod_i G_S) / \prod_i Aut(G_{T_i})$
	\end{center} 
	
	\textit{where $Aut_{b.a.}(\prod_i G_S)$ is the monoid of automorphisms of the form $\prod_i \alpha_i$ which are $S$-bounded above, see Definition 1.4.1.}

We will proceed, in Section 1.2, by reformulating the double coset result of $[1, \text{Prop. 7.5.2}]$ in terms of the 'better behaved' fracture square with base $\omega$ instead of $\Delta$. We will call this fracture square the 'formal' fracture square. This leads naturally to the notion of the 'restricted genus' of $G$, which turns out to be equivalent to the genus of $G$, and a double coset formula for the restricted genus of $G$.

To finish this introductory section we record two lemmas on nilpotent groups that will help us on our way in the remainder of these notes. The first is a consequence of the fact that the group $G$ we are working with is finitely $T$-generated:

\begin{Lemma} Let $G$ be an $f\mathbb{Z}_T$-nilpotent group with reference diagram as above. Then: \\
	i) $G$ is $T$-Noetherian; that is $G$ satisfies the ascending chain condition for $T$-local subgroups, \\
	ii) $\tilde{\pi}$ is a monomorphism, \\
	iii) $G_{T_i}$ has no $(T_i - S)$-torsion for all but finitely many $i$. Equivalently, $\phi_i$ is a monomorphism for all but finitely many $i$.
\end{Lemma}

\textbf{Proof:} i) This follows in the abelian case from the fact that $\mathbb{Z}_T$ is Noetherian and in the nilpotent case by induction up a suitable central series for $G$, \\
ii) It suffices to prove that $\prod \phi_i: \prod G_{T_i} \to \prod G_S$ is an $S$-monomorphism. This will follow from iii) and the fact that each $\phi_i$ is an $S$-monomorphism, \\
iii) Let $P = \{p_1,...,p_k\}$ be a finite set of prime numbers and define:

\begin{center}
	$G^{P} = \{g \in G | \ g^p = 1 \ \text{for some product p of primes in P}\}$
\end{center}

Then $G^P$ is a $T$-local subgroup of $G$, for example by Lemma 1.1.2 below. Since $G$ is $T$-Noetherian it follows that there is a finite set of primes $Q$ such that if $g^n = 1$ for some $n \in \mathbb{N}$, then $g^q = 1$ for some product of primes in $Q$. Now suppose that $T_i$ does not contain any primes in $Q$. If $a \in G_{T_i}$ is such that $a^s = 1$ for some product of primes in $(T_i - S)$, let $t_1$ be a $T_i$-number such that $a^{t_1} = \psi_i(g)$ for some $g \in G$. We have that $\psi_i(g^s) = 1$ and so there is a $T_i$-number $t_2$ such that $g^{st_2} = 1$. Since $s$ is coprime to each of the primes in $Q$, it follows that $g^{t_2} = 1$ and, therefore, that $a^{t_1t_2} = 1$. Since $G_{T_i}$ is $T_i$-local, it follows that $a = 1$, as desired.  $\strut\hfill \square$ 

The next lemma is a result about nilpotent groups which is essentially the same as $[2, \text{Thm.} \ 3.25]$.

\begin{Lemma} Let $G$ be a nilpotent group of nilpotency class $c$, $H$ a subgroup of $G$, and $A$ a set of elements of $G$ such that there exists an $s \in \mathbb{N}$ such that $a \in A \implies a^s \in H$. Then, if $g \in G$ is in the image of the free group generated by elements of $A$ and $H$, $g^{s^d} \in H$ where $d = \frac{1}{2} c(c+1)$.
	
\end{Lemma}

\textbf{Proof:} Let $K$ be the image in $G$ of the free subgroup generated by elements of $A$ and $H$. Then $K$ has nilpotency class $e \leq c$. Let:

\begin{center}
	$1 = \Gamma^e K \subset ... \subset \Gamma^1 K \subset \Gamma^0 K = K$
\end{center}

be the lower central series of $K$. Suppose that $k \in K$ is of the form $xh$ where $x \in \Gamma^i K$ and $h \in H$. Recall that $\frac{\Gamma^i K}{\Gamma^{i+1} K}$ is an abelian group generated by commutators of the form $[z_0,....,z_i]$ where each $z_i \in A \cup H$. Since the commutators are bilinear and $\frac{\Gamma^i K}{\Gamma^{i+1} K}$ is a central subgroup of $\frac{K}{\Gamma^{i+1}K}$, it follows that $k^{s^{i+1}} = yh^{'}$ for some $y \in \Gamma^{i+1} K$, $h^{'} \in H$. $\hfill \square$

\newpage

\subsection{The restricted genus of G}	
	
Recall the reference diagram from the beginning of the introductory section. To define the restricted genus, we start by observing that the isomorphism class of $\omega$ is completely determined by the group $G$ and does not depend on the choices of localisations $\psi_i, \phi$ and $\sigma$. For clarity, here we mean $\omega^{'}$ is isomorphic to $\omega$ if there exists a commutative square with horizontal maps $\omega,\omega^{'}$ and vertical maps isomorphisms. However, as we proceed we will need to restrict attention to a certain subgroup of automorphisms of $(\prod G_{T_i})_S$, which we will call 'diagonal' automorphisms. These are analogous to the automorphisms of $\prod(G_S)$ of the form $\prod_i \alpha_i$:

\begin{Definition} $DAut((\prod G_{T_i})_S)$ is the subgroup of $Aut((\prod G_{T_i})_S)$ consisting of automorphisms $\alpha$ such that, for every $j \in I$, under the canonical identification of $(\prod G_{T_i})_S$ with $G_{T_j} \times (\prod_{i \neq j} G_{T_i})_S$ determined by our reference diagram and any localisation of the form $\phi_j \times \phi^{'}$ of $G_{T_j} \times (\prod_{i \neq j} G_{T_i})$, $\alpha = \alpha_j \times \beta$ for some automorphisms $\alpha_j$ of $G_S$ and $\beta$ of $(\prod_{i \neq j} G_{T_i})_S$. \end{Definition}

Note that if $\alpha$ is a diagonal automorphism of $(\prod G_{T_i})_S$, then there is a commutative diagram:

\[ \begin{tikzcd}
(\prod G_{T_i})_S \arrow{r}{\tilde{\pi}} \arrow[swap]{d}{\alpha} & \prod G_S  \arrow{d}{\prod \alpha_i} \\
(\prod G_{T_i})_S \arrow[swap]{r}{\tilde{\pi}} & \prod G_S \\
\end{tikzcd}
\]  

Since $\tilde{\pi}$ is a monomorphism, it follows that there is an injective homomorphism $DAut((\prod G_{T_i})_S) \to \prod Aut(G_S)$. We will now show that the image of this map is the subgroup of '$S$-bounded' automorphisms of $\prod G_S$ and that all automorphisms of $\prod G_S$ induced by an automorphism of $G_S$ by the diagonal map are $S$-bounded. It follows that $Aut(G_S)$ defines a subgroup of $DAut((\prod G_{T_i})_S)$. 

\begin{Definition} An automorphism $\alpha = \prod \alpha_i \in \prod Aut(G_S)$ is said to be $S$-bounded if there exists an $S$-number $s$ such that for all $(g_i) \in \prod G_{T_i}$, $\alpha \circ (\prod \phi_i) (g_i^s) \in im(\prod(\phi_i))$ and $\alpha^{-1} \circ (\prod \phi_i) (g_i^s) \in im(\prod(\phi_i))$. \end{Definition}

\begin{Lemma} An automorphism $\alpha \in \prod Aut(G_S)$ is the image of a diagonal automorphism $\beta$ iff $\alpha$ is $S$-bounded. \end{Lemma}

\textbf{Proof:} ($\implies$) First suppose that $\alpha$ is the image of a diagonal automorphism $\beta$. Let $A$ be a finite set of $T$-generators for $G$. Since $\phi$ is an $S$-epimorphism, there exists an $S$-number $s$ such that for all $a \in A$,  $\beta \omega \sigma (a^s)$ and $\beta^{-1} \omega \sigma (a^s) \in im(\phi)$. Now $\psi_i(A)$ is a finite set of $T_i$ generators for $G_{T_i}$ and $im(\phi_i)$ is a $T_i$-local subgroup of $G_S$. It follows, by Lemma 1.1.2, that if $g_i \in G_{T_i}$ then $\alpha_i \phi_i (g_i^{s^d}), \alpha_i^{-1} \phi_i (g_i^{s^d}) \in im(\phi_i)$, where $d = \frac{1}{2} c(c+1)$ for $c$ the nilpotency class of $G$. Since $d$ is independent of $i$, it follows that $\alpha$ is $S$-bounded as desired.

($\impliedby$) Suppose that $\alpha$ is $S$-bounded. Let $F \subset I$ be the finite subset of $I$ such that $\phi_i$ is not a monomorphism for $i \in F$. If $i \notin F$, let $H_i$ be the subgroup of $G_{T_i}$ consisting of $g_i$ such that $\alpha_i \phi_i (g_i) \in im(\phi_i)$. Define unique homomorphisms $f_i$ such that the following square commutes:

\[ \begin{tikzcd}
\prod_{i \notin F} H_i \arrow{r}{\prod f_i} \arrow[swap]{d}{(\prod \phi_i) \circ \iota} & \prod_{i \notin F} G_{T_i} \arrow{d}{\prod \phi_i} \\
\prod_{i \notin F} G_S \arrow[swap]{r}{\prod \alpha_i} & \prod_{i \notin F} G_S \\
\end{tikzcd}
\]  

where $\iota$ is the inclusion of $\prod H_i$ into $\prod G_{T_i}$. Since $\alpha$ is $S$-bounded, $\iota$ is an $S$-isomorphism. Similarly, the image of the monomorphism $f := \prod f_i$ is $\{(g_i) \ | \ \forall i \  \alpha_i^{-1} \phi_i(g_i) \in im(\phi_i)\}$ and so $f$ is also an $S$-isomorphism. Let $\phi_{I/F}$ be a localisation of $\prod_{i \notin F} G_{T_i}$ so there is an induced isomorphism $(\prod G_{T_i})_S \cong \prod_{i \in F} G_S \times (\prod_{i \notin F} G_{T_i})_S$ coherent with our reference diagram. Since the vertical arrows in the diagram below are localisations at $S$, we can define the localisation $f_S$ as:

\[ \begin{tikzcd}
\prod_{i \notin F} H_i \arrow{r}{f} \arrow[swap]{d}{\phi_{I/F} \iota} & \prod_{i \notin F} G_{T_i} \arrow{d}{\phi_{I/F}} \\
(\prod_{i \notin F} G_{T_i})_S \arrow[swap]{r}{f_S} & (\prod_{i \notin F} G_{T_i})_S \\
\end{tikzcd}
\]  

Since $f$ is an $S$-isomorphism, $(\prod \alpha_i) \times f_S$ defines an automorphism of $\prod_{i \in F} G_S \times (\prod_{i \notin F} G_{T_i})_S$. Uniqueness of localisations now implies that if we define $\beta \in Aut((\prod G_{T_i})_S)$ to correspond to $(\prod \alpha_i) \times f_S$ under the induced isomorphism given above, then $\beta$ is a diagonal automorphism whose image is $\alpha$. $\hfill \square$

\begin{Lemma} If $\alpha \in Aut(G_S)$, then $\prod \alpha$ is $S$-bounded. \end{Lemma}

\textbf{Proof:} Let $A$ be a finite set of $T$-generators for $G$. Since $\sigma$ is an S-localisation, there exists an $S$-number $s$ such that for all $a \in A$, $\alpha \sigma (a^s), \alpha^{-1} \sigma (a^s) \in im(\sigma)$. Since, for all $i$, $\sigma = \phi_i \psi_i$ and $\psi_i(A)$ is a finite set of $T_i$ generators for $G_{T_i}$, this implies, by Lemma 1.1.2, that for all $g_i \in G_{T_i}$, $\alpha \phi_i(g_i^{s^d})$ and $\alpha^{-1} \phi_i(g_i^{s^d}) \in im(\phi_i)$, where $d = \frac{1}{2}c(c+1)$ is independent of $i$. It follows that $\prod \alpha$ is $S$-bounded. $\hfill \square$

We are now in a position to prove a double coset result for the formal fracture square. However, not all elements of the extended genus will 'factor through' the formal fracture square. This motivates the following definition:

\begin{Definition} The restricted genus of $G$ is defined to be the set of isomorphism classes of T-local nilpotent groups $H$ such that $H_{T_i} \cong G_{T_i}$, and if $(\varphi_i): H \to \prod G_{T_i}$ is a map such that each $\varphi_i$ is a localisation at $T_i$, then the localisation $\omega^{'}$ in the diagram below:

\[ \begin{tikzcd}
H \arrow{r}{(\varphi_i)} \arrow[swap]{d}{\mu} & \prod G_{T_i}  \arrow{d}{\phi} \\
G_S \arrow[swap]{r}{\omega^{'}} & (\prod G_{T_i})_S \\
\end{tikzcd}
\]  

is diagonally isomorphic to $\omega$ in the sense that there is a commutative square with horizontal maps $\omega, \omega^{'}$, the left vertical map an automorphism of $G_S$ and the right vertical map a diagonal automorphism of $(\prod G_{T_i})_S$. \end{Definition}

Note that this definition is independent of the choice of localisation of $H$ at $S$, $\mu$. This definition will turn out to be a temporary one as we will show the restricted genus of $G$ is equivalent to the genus of $G$ in Section 1.3. 

The proof of the following theorem is essentially the same as the intended proof in $[1, \text{Prop. 7.5.2}]$:

\begin{Theorem}There is a 1-1 correspondence between the restricted genus of $G$ and the double coset:
	
	\begin{center}
		$Aut(G_S)$ \textbackslash \ $DAut((\prod G_{T_i})_S) \ / \ \prod Aut(G_{T_i})$	
	\end{center}
	
	This correspondence sends a diagonal automorphism $\alpha$ to the pullback group of $\alpha \phi$ along $\omega$. \end{Theorem}

\textbf{Proof:} Let $\alpha$ be a diagonal automorphism. Our first task is to show that the pullback of $\alpha \phi$ along $\omega$ is in the restricted genus of $G$. The stated pullback is isomorphic to the following pullback:

\[ \begin{tikzcd}
P \arrow{r}{(\varphi_i)} \arrow[swap]{d}{\mu} & G_{T_j} \times (\prod_{i \neq j} G_{T_i}) \arrow{d}{(\alpha_j \phi_j) \times (\alpha^{'} \phi^{'})} \\
G_S \arrow[swap]{r}{1 \times \omega^{'}} & G_S \times (\prod_{i \neq j} G_{T_i})_S \\
\end{tikzcd}
\]  

where $\phi^{'}$ is some localisation of $\prod_{i \neq j} G_{T_i}$ and the rest of the diagram is filled in to be coherent with our reference diagram and we use the fact that $\alpha$ is a diagonal automorphism. Now $1 \times (\alpha^{'} \phi^{'})$ is a localisation at $T_j$, so we can localise the pullback with respect to this localisation, the identity localisations of the bottom row, and an arbitrary $T_j$-localisation of $P$, to obtain a pullback:

\[ \begin{tikzcd}[column sep=7em]
P_{T_j} \arrow{r}{(\varphi_j)_{T_j} \times ((\varphi_i))_{T_j}} \arrow[swap]{d}{\mu_{T_j}} & G_{T_j} \times (\prod_{i \neq j} G_{T_i})_S \arrow{d}{(\alpha_j \phi_j) \times 1} \\
G_S \arrow[swap]{r}{1 \times \omega^{'}} & G_S \times (\prod_{i \neq j} G_{T_i})_S \\
\end{tikzcd}
\]  

So $P_{T_j} \cong G_{T_j}$, the top horizontal map is isomorphic to $1 \times (\omega^{'}\alpha_j\phi_j)$ and the left hand vertical map is isomorphic to $\alpha_j \phi_j$. In particular, $\varphi_j$ is a $T_j$-localisation for all $j$ and $\mu$ is an $S$-localisation. Moreover, since $\alpha$ is a diagonal automorphism, the localisation of $(\varphi_i)$ at $S$ is diagonally isomorphic to $\omega$ and so $P$ represents an element of the restricted genus of $G$, as desired. 

It is clear that the isomorphism class of $P$ is the same for any element $\beta$ in the same double coset equivalence class as $\alpha$. Therefore, we have a well defined map from the double coset to the restricted genus and it remains to show injectivity and surjectivity. Surjectivity follows from the definition of the restricted genus. For injectivity, suppose that $\alpha, \beta$ are diagonal automorphisms and we have pullbacks:

\[ \begin{tikzcd}
P \arrow{r}{(\varphi_i)} \arrow[swap]{d}{\mu} & \prod G_{T_i} \arrow{d}{\alpha \phi} & P \arrow{r}{(\bar{\varphi}_i)} \arrow[swap]{d}{\bar{\mu}} & \prod G_{T_i} \arrow{d}{\beta \phi} \\
G_S \arrow[swap]{r}{\omega} & (\prod G_{T_i})_S & G_S \arrow[swap]{r}{\omega} & (\prod G_{T_i})_S \\
\end{tikzcd}
\]  

By uniqueness of localisations, there is an automorphism $\gamma$ of $\prod G_{T_i}$ such that $\gamma (\varphi_i) = (\bar{\varphi}_i)$. Therefore, since we only care about equivalence classes in the double coset we may assume that $\varphi_i = \bar{\varphi}_i$ for all $i$. Similarly, there is an automorphism $\gamma^{'}$ of $G_S$ such that $\gamma^{'} \mu = \bar{\mu}$ and so we can reduce to the case $\mu = \bar{\mu}$. Now, for all $i$, $\alpha_i \phi_i$ and $ \beta_i \phi_i$ are both the unique factorisation of $\mu$ through $\varphi_i$. Therefore, $(\tilde{\pi} \alpha) \phi = (\tilde{\pi} \beta) \phi$. By uniqueness of factorisation through $\phi$, we have $\tilde{\pi} \alpha = \tilde{\pi} \beta$ and, since $\tilde{\pi}$ is a monomorphism, this implies that $\alpha = \beta$. $\hfill \square$

Finally, for this section we relate our result to the alternative fracture square with $\Delta$ along the bottom instead of $\omega$.

\begin{Corollary} The restricted genus of $G$ is in 1-1 correspondence with the double coset:
	
	\begin{center}
		$Aut(G_S) \text{\textbackslash} Aut_{b}(\prod_i G_S) / \prod_i Aut(G_{T_i})$
	\end{center} 
	
	where $Aut_{b}(\prod_i G_S)$ is the subgroup of automorphisms of the form $\prod_i \alpha_i$ which are $S$-bounded.
\end{Corollary}

\newpage

\subsection{Restricted genus equals the genus}

We begin with the following lemma, which is essentially the same as $[1, \text{Prop. 7.4.3}]$. 

\begin{Lemma} If $I$ is a finite indexing set and $H$ is a $T$-local nilpotent group such that $H_{T_i}$ is finitely $T_i$-generated for all $i \in I$, then $H$ is finitely $T$-generated.
	
\end{Lemma}

\textbf{Proof:} Let $H_0 \subset H_1 \subset ...$ be an ascending chain of $T$-local subgroups of $H$. For each $i$, let $\psi_i$ denote a $T_i$-localisation of $H$ and let $H_j^{i}$ denote the $T_i$-local subgroup of $H_{T_i}$ generated by $\psi_i(H_j)$. Choose an integer $N$ such that $H_0^i \subset H_1^i \subset ...$ terminates at $H_N^i$ for all $i$. Now let $n \geq N$; we claim that $H_n = H_N$. If $h \in H_n$, then there exists a $T_i$-number $t_i$ and a $k \in H_N$ such that $\psi_i(h^{t_i}) = \psi_i(k)$. It follows that there is a $T_i$-number $s_i$ such that $(h^{t_i}k^{-1})^{s_i} = 1$. Since the set of $g \in H$ such that there exists a $T_i$-number $s$ such that $g^s \in H_N$ is a subgroup of $G$ which contains $h^{t_i}k^{-1}$ and $k$, it follows that there is a $T_i$-number $r_i$ such that $h^{r_i} \in H_N$. Now any common factor of each of the $r_i$ lies outside of $T$ and so there is a $T$-number $r$ such that $h^r \in H_N$. Since $H_N$ is $T$-local, it follows that $h \in H_N$ as desired. So $H$ is $T$-Noetherian $\implies$ $H$ is finitely $T$-generated. $\strut\hfill \square$

Now suppose that $H$ is in the restricted genus of $G$. We again consider the finite subset of $I$ consisting of $i$ such that $\phi_i$ is not a monomorphism. Then $H$ fits into a diagram of the form:

\[ \begin{tikzcd}
H \arrow{r} \arrow{d} & (\prod_{i \notin F} G_{T_i}) \times (\prod_{j \in F} G_{T_j}) \arrow{d}{1 \times (\prod ( \alpha_j \phi_j))} \\
P \arrow{r} \arrow{d} & (\prod_{i \notin F} G_{T_i}) \times (\prod_{j \in F} G_S) \arrow{r} \arrow{d}{(\prod (\alpha_i \phi_i)) \times 1} & \prod_{i \notin F} G_{T_i} \arrow{d} \\
G_S \arrow[swap]{r}{\Delta \times \Delta} & (\prod_{i \notin F} G_{S}) \times (\prod_{j \in F} G_S) \arrow{r} & \prod_{i \notin F} G_{S}\\
\end{tikzcd}
 \]

where $\alpha$ is $S$-bounded and each of the squares is a pullback. Consider the localisation of the diagram at $T^{'} = \cup_{i \notin F} T_i$. The objects in the bottom two rows are all $T^{'}$-local. If $j \in F$, $T_j \cap T^{'} = S$, so $\alpha_j \phi_j$ is a $T^{'}$-localisation. It follows that $P$ is a $T^{'}$-localisation of $H$. In light of Lemma 1.3.1, if we want to show that $H$ is finitely $T$-generated it suffices to show that $P$ is finitely $T^{'}$-generated. Note also that $P$ is in the restricted genus of $G_{T^{'}}$. In this way we can reduce the following lemma to the case where $\phi_i$ is a monomorphism for all $i$. 

\begin{Lemma} If $H$ is in the restricted genus of $G$, then $H$ is finitely $T$-generated.
	
\end{Lemma}

\textbf{Proof:} As discussed above, we can reduce to the case where $\phi_i$ is a monomorphism for all $i$. By Corollary 1.2.7, there is an $S$-bounded automorphism $\alpha = \prod \alpha_i$ of $\prod G_S$ such that we have a pullback square:

\[ \begin{tikzcd}
H \arrow{r} \arrow{d} & \prod G_{T_i} \arrow{d}{\prod \alpha_i \phi_i} \\
G_S \arrow[swap]{r}{\Delta} & \prod G_S \\
\end{tikzcd}
\]

Let $K$ be the $T$-local subgroup of $H$ consisting of pairs $(x,(g_i))$ with $x \in G_S, g_i \in G_{T_i}$ such that, for all $i$, $\alpha_i \phi_i(g_i) = x$ and $x \in im(\phi_i)$, say $x = \phi_i(a_i)$. Then there is an injective group homomorphism $K \to G$ sending $(x,(g_i))$ to $(x,(a_i))$. Since $G$ is finitely $T$-generated so is $K$, and since $\alpha$ is $S$-bounded there exists an $S$-number $s$ such that if $h \in H$, then $h^s \in K$. Consider a $T$-subnormal series for $K$:

\begin{center}
	$K = K_0 \subset K_1 \subset ... \subset K_m = H$
\end{center}

If we localise at $T_i$, then all of the groups in the chain become finitely $T_i$-generated. Moreover, $(\frac{K_{j+1}}{K_j})_{T_i}$ is a finitely $T_i$-generated nilpotent group such that if $k \in (\frac{K_{j+1}}{K_j})_{T_i}$, then $k^s = 1$. For all but finitely many $i$ this implies that $(\frac{K_{j+1}}{K_j})_{T_i}$ is trivial. For the remaining $i$, $(\frac{K_{j+1}}{K_j})_{T_i}$ is finitely $T_i$-generated (in fact it is finite). Therefore, using a global to local fracture square we see that $\frac{K_{j+1}}{K_j}$ is finitely $T$-generated (in fact it is finite). Inductively, it follows that $H$ is finitely $T$-generated (and $K$ is a subgroup of finite index in $H$). $\hfill \square$

It remains to prove that every element of the genus is in the restricted genus. We start with the following observation:

\begin{Lemma} If $H$ is in the same genus as $G$, then there is a finite subset $F$ of $I$ such that if $T^{'} = \cup_{i \notin F} T_i$, then $G_{T^{'}} \cong H_{T^{'}}$.
\end{Lemma}

\textbf{Proof:} It was observed in $[3, \text{Thm I.3.3}]$, that, since $G_S \cong H_S$, there is a finitely $T$-generated nilpotent group $P$ equipped with $S$-isomorphisms $f: P \to G$ and $g: P \to H$. In fact, we just need to consider the pullback:

\[\begin{tikzcd}
P \arrow{r} \arrow{d} & G \arrow{d}{\phi_{S}} \\
H \arrow[swap]{r}{\phi_S} & G_S \\
\end{tikzcd}
\]

to get the desired maps, where $\phi_S$ denotes a localisation at $S$. Since $G,H$ and $P$ are finitely $T$-generated, we can use Lemma 1.1.2 to show that there exists an $S$-number $s$ such that if $p \in ker(f)$ or $p \in ker(g)$, then $p^s = 1$ and, if $g \in G, h \in H$, then $g^s \in im(f), h^s \in im(g)$. This implies that if we take $T^{'}$ to be the union of the $T_i$ which don't contain any prime factors of $s$, then both $f$ and $g$ are $T^{'}$-isomorphisms, which implies the result. $\hfill \square$

\newpage

\begin{Lemma} If $H$ is in the genus of $G$, then $H$ is in the restricted genus of $G$.\end{Lemma}
	
\textbf{Proof:} Let $F$ be a finite subset of $I$ such that if $T' = \cup_{i \notin F} T_i$, then $H_{T^{'}} \cong G_{T^{'}}$. Let $\mu: H \to G_{T^{'}}$ and $\epsilon: G \to G_{T^{'}}$ be $T^{'}$-localisations. Then there are unique factorisations of $\sigma$ and $\psi_i$, for $i \notin F$, through $\epsilon$; denote them by $\sigma^{'}, \psi_i^{'}$. Since $H$ is finitely $T$-generated, we can form a global to local fracture square of the form:

\[ \begin{tikzcd}[column sep = 6em]
H \arrow{r}{(\psi_i^{'} \mu) \times (\varphi_j)} \arrow[swap]{d}{\sigma^{'} \mu} & (\prod_{i \notin F} G_{T_i}) \times (\prod_{j \in F} G_{T_j}) \arrow{d}{(\phi_i) \times (\alpha_j \phi_j)} \\
G_S \arrow[swap]{r}{\Delta \times \Delta} & (\prod_{i \notin F} G_S) \times (\prod_{j \in F} G_S) \\
\end{tikzcd}
\]

where $\varphi_j$ is any $T_j$-localisation of $H$ and $\alpha_j \in Aut(G_S)$. Since $F$ is finite, $1 \times (\alpha_j)$ is $S$-bounded, so $H$ is in the restricted genus by Corollary 1.2.7. $\hfill \square $

\newpage

\subsection{Relation to the extended genus}

The aim of this section is to show that if $\alpha = \prod \alpha_i$ is an automorphism of $\prod G_S$, then in the pullback diagram below:

\[ \begin{tikzcd}
H \arrow{r}{(\varphi_i)} \arrow[swap]{d}{\mu} & \prod G_{T_i} \arrow{d}{\prod \alpha_i \phi_i} \\
G_S \arrow[swap]{r}{\Delta} & \prod G_S \\
\end{tikzcd}
\]

$\varphi_i$ is a $T_i$-localisation for all $i$ iff $\alpha$ is $S$-bounded above in the following sense:

\begin{Definition} An automorphism $\alpha = \prod \alpha_i \in \prod Aut(G_S)$ is said to be $S$-bounded above if there exists an $S$-number $s$ such that for all $i$ and for all $g_i \in G_{T_i}$, $\alpha_i^{-1} \phi_i(g_i^s) \in im(\phi_i)$. \end{Definition}

We start with:

\begin{Lemma} If $\varphi_i$ is a $T_i$-localisation for all $i$, then $\alpha$ is $S$-bounded above.
\end{Lemma}

\textbf{Proof:} Let $A$ be a finite $T$-generating set for $G$. Since all the $\varphi_i$ are $T_i$-localisations, $\mu$ is an $S$-localisation. It follows that there exists an $S$-number $s$ such that for all $a \in A$, $\sigma(a^s) \in im(\mu) \subset im(\alpha_i \phi_i)$ for all $i$. Since $\psi_i(A)$ is a finite $T_i$-generating set for $G_{T_i}$, it follows from Lemma 1.1.2 that if $g_i \in G_{T_i}$, then $\phi_i (g_i^{s^d}) \in im(\alpha_i \phi_i)$, where $d = \frac{1}{2} c(c+1)$, for $c$ the nilpotency class of $G$. $\hfill \square$

For the reverse direction, we start with the following observation which does not require $\alpha$ to be $S$-bounded above:

\begin{Lemma} $\mu$ is an $S$-monomorphism.	
\end{Lemma}

\textbf{Proof:} Since $\prod \phi_i$ is an $S$-monomorphism, so is $\prod \alpha_i \phi_i$, and the pullback of an $S$-monomorphism is an $S$-monomorphism. $\hfill \square$

\begin{Lemma} If $\alpha$ is $S$-bounded above, then $\mu$ is an $S$-epimorphism, hence an $S$-localisation.
\end{Lemma}

\textbf{Proof:} If $x \in G_S$, then since $\sigma$ is an $S$-localisation, there exists an $S$-number $r$ such that $x^r \in im(\sigma) \subset im(\phi_i)$ for all $i$. Since $\alpha$ is $S$-bounded, there exists an $S$-number $s$ such that $x^{rs} \in im(\alpha_i \phi_i)$ for all $i$. It follows that $x^{rs}$ is in the image of $\mu$ by the definition of a pullback. $\hfill \square$

We now have:

\begin{Lemma} If $\alpha$ is $S$-bounded above, then $\varphi_i$ is a $T_i$-localisation for all $i$.	
\end{Lemma}

\textbf{Proof:} If $h \in H$ and $\varphi_i(h) = 1$, then $\mu(h) = 1$ and so there exists an $S$-number $s$ such that $h^s = 1$. Write $s$ as a product of a $T_i$-number $t$ and a product of primes in $T_i$, $r$. Then if $j \neq i$, $\varphi_j(h^t) = 1$, since $G_{T_j}$ is $T_j$-local and $T_i \cap T_j = S$. Clearly $\varphi_i(h^t) = 1$, so it follows that $\varphi_i$ is a $T_i$-monomorphism.

Now suppose that $g_i \in G_{T_i}$ and let $x = \alpha_i \phi_i (g_i)$. Since $\mu$ is an $S$-localisation, there exists an $S$-number $s$ and $h \in H$ such that $x^s = \mu(h)$. Write $s$ as a product of a $T_i$-number $t$ and a product of primes in $T_i$, $r$. If $j \neq i$, then the image of $\alpha_j \phi_j$ is a $T_j$-local subgroup of $G_S$ and so $x^t \in im(\alpha_j \phi_j)$ for all $j \neq i$. Since $x^t$ is also in $im(\alpha_i \phi_i)$, it follows that $g_i^{t}$ is in the image of $\varphi_i$ by the definition of a pullback. $\hfill \square$

It  is now possible to prove the following double coset formula, by repeating the arguments of $[1, \text{Prop. 7.5.2}]$ or Theorem 1.2.6 above:

\begin{Theorem} The extended genus of $G$ is in 1-1 correspondence with the double coset:
	
	\begin{center}
		$Aut(G_S) \text{\textbackslash} Aut_{b.a.}(\prod_i G_S) / \prod_i Aut(G_{T_i})$
	\end{center} 

where $Aut_{b.a.}(\prod_i G_S)$ is the monoid of automorphisms of the form $\prod_i \alpha_i$ which are $S$-bounded above.
\end{Theorem}

\newpage

\section{References}

$[1]$ -  J. P. May and K. Ponto. \textit{More concise Algebraic Topology: localization, completion, and model categories}. Chicago lectures in mathematics. University of Chicago Press, Chicago, Ill. 2012 

$[2]$ -  R. B.Warfield, Jr. \textit{Nilpotent groups}, volume 513 of Lecture notes in mathematics. SpringerVerlag, Berlin, 1976.

$[3]$ - P. Hilton, G. Mislin, and J. Roitberg. \textit{Localization of nilpotent groups and spaces}. Volume 15 in North-Holland mathematics studies. North-Holland Publishing, Amsterdam, 1975.

\end{document}